\documentclass [11pt,twoside,a4paper]{article}
\usepackage{amsfonts}
\usepackage{amsthm}
\usepackage{amsmath}
\usepackage{amstext}
\usepackage{amssymb}
\usepackage{mathrsfs}
\usepackage{amscd}
\usepackage{xypic}
\usepackage{epsf}              
\usepackage{graphicx}          
\usepackage{fancybox}          
\usepackage{color}             
\usepackage{fancyhdr}
\usepackage[hang,footnotesize]{caption2}  

\numberwithin{equation}{section}

\def\mathcal{\mathscr}
\newfont{\aaa}{cmb10 at 19pt}
\newfont{\bbb}{cmb10 at 11pt}

\def\v1{\vspace{1mm}}

\def\leq{\leqslant}

\def\geq{\geqslant}

\def\dfrac{\displaystyle\frac}

\newcommand{\beq}{\begin{equation}}
\newcommand{\eeq}{\end{equation}}
\newcommand{\bey}{\begin{eqnarray}}
\newcommand{\eey}{\end{eqnarray}}
\newcommand{\beyy}{\begin{eqnarray*}}
\newcommand{\eeyy}{\end{eqnarray*}}

\makeatletter
\def\@evenhead{
\vbox{\hbox to \textwidth {}{\hspace{0mm}{\footnotesize
\thepage}}{\hspace{8cm} {\footnotesize {Manli Song and Jiale Yang}}}
\protect\vspace{1truemm}\relax \hrule depth0pt
height0.15truemm width\textwidth}}
\def\@evenfoot{}
\def\@oddhead{\vbox{\hbox to \textwidth
{{\hspace{0cm}{\footnotesize Decay estimates for a class of wave equations on the Heisenberg group}\hfill{\footnotesize
\thepage}}\hspace{0mm}}{} \protect\vspace{1truemm}\relax\hrule
depth0pt height0.15truemm width\textwidth}}
\def\@oddfoot{}
\makeatother

\begin{document}



\setcounter{page}{1}
\qquad\\[8mm]

\noindent{\aaa{Decay estimates for a class of wave equations on the Heisenberg group}}\\[1mm]

\noindent{\bbb Manli Song and Jiale Yang}\\[-1mm]

\noindent\footnotesize{School of Mathematics and Statistics, Northwestern Polytechnical University, Xi'an, Shaanxi 710129, China}\\[6mm]

\normalsize\noindent{\bbb Abstract}\quad In this paper, we study a class of dispersive wave equations on the Heisenberg group $H^n$.  Based on the group Fourier transform on $H^n$, the properties of the Laguerre functions and the stationary phase lemma, we establish the decay estimates for a class of dispersive semigroup on $H^n$ given by $e^{it\phi(\mathcal{L})}$, where $\phi: \mathbb{R}^+ \to \mathbb{R}$ is smooth, and $\mathcal{L}$ is the sub-Laplacian on $H^n$. Finally, using the duality arguments, we apply the obtained results to derive the Strichartz inequalities for the solutions of some specific equations, such as the fractional Schr\"{o}dinger equation, the fractional wave equation and the fourth-order Schr\"{o}dinger equation.
\vspace{0.3cm}

\footnotetext{
Corresponding author: Jiale Yang, E-mail:
yangjiale2210@163.com}

\noindent{\bbb Keywords}\quad  Decay estimates, Heisenberg group, Sub-Laplacian.\\
{\bbb MSC}\quad {22E25, 33C45, 35H20, 35B40}\\[0.4cm]

\newtheorem{theorem}{Theorem}[section]
\newtheorem{preliminaries}{Preliminaries}[section]
\newtheorem{definition}{Definition}[section]
\newtheorem{main result}{Main Result}[section]
\newtheorem{lemma}{Lemma}[section]
\newtheorem{proposition}{Proposition}[section]
\newtheorem{corollary}{Corollary}[section]
\newtheorem{remark}{Remark}[section]


\section{Introduction}

The aim of this paper is to investigate the decay estimates for a class of dispersive equations on the Heisenberg group $H^n$
	\begin{equation}\label{HSEquation}
	\begin{cases}
	i\partial_tu+\phi(\mathcal{L}) u=f,\\
	u|_{t=0}=u_0,
	\end{cases}
	\end{equation}
where $\phi: \mathbb{R}^+ \to \mathbb{R}$ is smooth, and $\mathcal{L}$ is the sub-Laplacian on the Heisenberg group $H^n$. 

In the past decades, Strichartz estimates have been very useful in the study of nonlinear partial differential equations. These estimates in the Euclidean setting have been proved for many dispersive equations, such as the wave equation and Schr\"{o}dinger equation (see \cite{GV, KT, Str}). To obtain Strichartz estimates, it involves basically two types of ingredients. The first one consists in estimating the decay in time on the evolution group associated with the free equation (i.e. with $f=0$). The second one consists of abstract arguments, which are mainly duality arguments. Therefore, the decay estimate plays a crucial role.

In 2008, Z. Guo, L. Peng and B. Wang\cite{GPW2008} used a unified way to study the decay for a class of dispersive semigroup on the Euclidean space $\mathbb{R}^n$ given by $e^{it\phi(\Delta)}$, where $\Delta=-\sum\limits_{j=1}^n\partial_{x_j}^2$ is the Laplacian on $\mathbb{R}^n$. Many dispersive wave equations reduce to this type, for instance, the Schr\"{o}dinger equation ($\phi(r)=r$), the wave equation ($\phi(r)=r^{\frac{1}{2}}$), the fractional Schr\"{o}dinger equation ($\phi(r)=r^\alpha$) ($0<\alpha<1$), the fourth-order Schr\"{o}dinger equation ($\phi(r)=r^2+r$), etc. When $\phi$ is a homogeneous function of order $m$, namely, $\phi(\lambda r)=\lambda^m\phi(r),\,\forall \lambda>0$, the dispersive estimate can be easily obtained by a theorem of Littman and dyadic decomposition. However, it becomes very complicated when $\phi$ is not homogeneous since the scaling constants can not be effectively separated from the time.  To overcome the difficulty, they applied frequency localization by separating $\phi$ between high and low frequency in different scales. We shall assume that $\phi: \mathbb{R}^+ \to \mathbb{R}$ is smooth and satisfies:  \\

(H1)~There exists $m_1>0$, such that for any $\alpha \geqslant 2$ and $\alpha \in \mathbb{N}$,
\begin{equation*}
|\phi'(r)| \sim r^{m_1-1}  , \quad |\phi^{(\alpha)}(r)| \lesssim r^{m_1-\alpha},\quad r \geqslant 1.
\end{equation*}

(H2)~There exists $m_2>0$, such that for any $\alpha \geqslant 2$ and $\alpha \in \mathbb{N}$,
\begin{equation*}
|\phi'(r)| \sim r^{m_2-1}  , \quad |\phi^{(\alpha)}(r)| \lesssim r^{m_2-\alpha},\quad 0<r<1.
\end{equation*}

(H3)~There exists $m_1>0$, such that
\begin{equation*}
|\phi''(r)| \sim r^{\alpha_1-2}, \quad r \geqslant1.
\end{equation*}

(H4)~There exists $m_2>0$, such that
\begin{equation*}
|\phi''(r)| \sim r^{\alpha_2-2}, \quad  0<r<1.
\end{equation*}

 (H1) and (H3) represent the homogeneous order of $\phi$ in high frequency. Similarly, the homogeneous order of $\phi$ in low frequency is described by (H2) and (H4).
\begin{remark} 
	 Throughout this paper, $A\lesssim B$ means that $A\leqslant CB$, and $A\sim B$ stands for $C_1B\leq A\leq C_2B$, where $C$, $C_1$, $C_2$ denote positive universal constants.
\end{remark}

Choose a Littlewood-Paley dyadic  partition of unity $\{\Phi_j\}_{j\in Z}$ on $\mathbb{R}^n$ and define the Littlewood-Paley projector by $\Delta_j=\mathcal{F}\Phi_j(\xi)\mathcal{F}^{-1}$, where $\mathcal{F}$ is the Fourier transform on $\mathbb{R}^n$. Z. Guo, L. Peng and B. Wang\cite{GPW2008} obtained the following decay estimates in time for the semigroup $e^{it\phi(\Delta)}$:

	(1)~For $j \geqslant 0$, $\phi$ satisfies (H1), then
	\begin{equation*}
	\left\|e^{it\phi(\sqrt{-\Delta})}\Delta_ju_0\right\|_{L^\infty(\mathbb{R}^n)}\lesssim |t|^{-\theta}2^{j(n-m_1\theta)}\|u_0\|_{L^1(\mathbb{R}^n)},\,0\leq\theta\leq\frac{n-1}{2};
	\end{equation*}
In addition, if $\phi$ satisfies (H3), then 
	\begin{equation*}
	\left\|e^{it\phi(\sqrt{-\Delta})}\Delta_ju_0\right\|_{L^\infty(\mathbb{R}^n)}\lesssim|t|^{-\frac{n-1+\theta}{2}}2^{j(n-\frac{m_1(n-1+\theta)}{2}-\frac{\theta(\alpha_1-m_1)}{2})}\|u_0\|_{L^1(\mathbb{R}^n)},\, 0 \leq \theta \leq 1.
	\end{equation*}
	
	(2)~For $j < 0$, $\phi$ satisfies (H2), then
		\begin{equation*}
		\left\|e^{it\phi(\sqrt{-\Delta})}\Delta_ju_0\right\|_{L^\infty(\mathbb{R}^n)}\lesssim |t|^{-\theta}2^{j(n-m_2\theta)}\|u_0\|_{L^1(\mathbb{R}^n)},\,0\leq\theta\leq\frac{n-1}{2};
		\end{equation*}
	In addition, if $\phi$ satisfies (H4), then 
\begin{equation*}
\left\|e^{it\phi(\sqrt{-\Delta})}\Delta_ju_0\right\|_{L^\infty(\mathbb{R}^n)}\lesssim|t|^{-\frac{n-1+\theta}{2}}2^{j(n-\frac{m_2(n-1+\theta)}{2}-\frac{\theta(\alpha_2-m_2)}{2})}\|u_0\|_{L^1(\mathbb{R}^n)},\, 0 \leq \theta \leq 1.
\end{equation*}

Besides, by using standard duality arguments, they applied the above estimates to some concrete  equations and established the corresponding Strichartz estimates. Their results could not only cover known results but also make some improvements and provide simple proof for some equations. 

Many authors are also interested in adapting the well known Strichartz estimates from the Euclidean setting to a more abstract setting, such as the Heisenberg type groups and metric measure spaces. In 2000, H. Bahouri, P. G\'{e}rard and C. J. Xu \cite{BGX2000} initially discussed the Strichartz estimates with the sub-Laplacian on the Heisenberg group $H^n$, by means of Besov spaces defined by a Littlewood-Paley decomposition related to the spectral of the sub-Laplacian. In their work, they showed such estimates existed for the wave equation, while surprisingly failed for the Schr\"{o}dinger equation.
To avoid the particular behavior of the Schr\"{o}dinger operator, one way is to replace the Heisenberg group by a bigger space, like the Heisenberg type group whose center dimension is bigger than $1$, on which it was shown in \cite{H2005} that  the Schr\"{o}dinger equation is dispersive. For more results on H-type groups and more generality of step $2$ stratified Lie groups, see \cite{BKG, Song, SZ}.  Another way is to adapt the Fourier transform restriction analysis initiated by R. S. Strichartz \cite{Str}  in the  Euclidean space. Using this approach, H. Bahouri, D. Barilari and I. Gallagher\cite{BBG2021} obtain an anisotropic Strichartz estimate for the Schr\"{o}dinger operator on $H^n$. For more results on the Heisenberg group, see \cite{FMV1, FV}. The case in the metric measure spaces is also considered, such as \cite{ BDXM2019, NR2005, R2008, S2013}.

Our purpose is to establish the decay estimates for a class of dispersive equations \eqref{HSEquation} on the Heisenberg group. The proof originates from the combination of H. Bahouri, P. G\'{e}rard and C. J. Xu \cite{BGX2000} and Z. Guo, L. Peng and B. Wang\cite{GPW2008}. When proving the dispersive estimates, the main difficulty is that the Fourier transform on the Heisenberg group is more complicated than in the Euclidean case. Our main results are in the following.

\begin{theorem}\label{ResultTime}
	Assume $\phi:\mathbb{R}^+ \to \mathbb{R}$ is smooth. $U_t=e^{it\phi(\mathcal{L})}$. $\varphi_j$ is Littlewood-Paley function introduced in section 2 and $N=2n+2$ is the homogeneous dimension of $H^n$. Then we have the following results.
	
	(1)~For $j \geq 0$, $\phi$ satisfies (H1), then
	\begin{equation}\label{res3-2}
	\|U_t \varphi_j\|_{L^\infty(H^n)} \leq C_n2^{Nj};
	\end{equation}
In addition, if $\phi$ satisfies (H3), then 
	\begin{equation}\label{res3-3}
	\|U_t \varphi_j\|_{L^\infty(H^n)} \leq C_n|t|^{-\frac{\theta}{2}}2^{j(N-\alpha_1\theta)},\,   0 \leq \theta \leq 1.
	\end{equation}
	
	(2)~For $j < 0$, $\phi$ satisfies (H2), then
	\begin{equation}
	\|U_t \varphi_j\|_{L^\infty(H^n)} \leq C_n2^{Nj};
	\end{equation}
In addition, if $\phi$ satisfies (H4), then 
	\begin{equation}
	\|U_t \varphi_j\|_{L^\infty(H^n)} \leq C_n|t|^{-\frac{\theta}{2}}2^{j(N-\alpha_2\theta)},\,
	0 \leq \theta \leq 1.
	\end{equation}	
\end{theorem}

This paper is organized as follows: In Section 2, we introduce some basic notations and definitions on the Heisenberg group. We also
recall the properties for spherical Fourier transform and the homogeneous Besov spaces on the Heisenberg group.  Section
3 is dedicated to describe the useful technical lemmas for proving the decay estimates and Strichartz estimates. In Section 4, we prove our main results Theorem \ref{ResultTime}. In Section 5, we exploit Theorem \ref{ResultTime} to obtain the Strichartz estimates for some concrete equations, such as the fractional Sch\"{o}dinger equation, the fractional wave equation and the fourth-order Schr\"{o}dinger equation. It is worth noticing that the list of applications is not exhaustive since we just intend to show the generality our approach.

\noindent\\[4mm]
\section{Prilimaries}
\subsection{Heisenberg groups}
In this section,we recall some fundamental conceptions and notations on the Heisenberg group. These facts are familiar and can be easily found in many references.

The Heisenberg group $H^n$ is a non-commutative nilpotent Lie group with the underlying mainfold $\mathbb{C}^n \times \mathbb{R}$. The element of $H^n$ is denoted by $(z,s)$ where $z=x+iy\in \mathbb{C}^n$ and $s\in \mathbb{R}$. The group law is given by 
	\begin{equation*}
	(z,s)\cdot (z',s')=(z+z',s+s'+2Im(z\cdot \bar{z'})),
	\end{equation*}
where $z\cdot \bar{z'}=z_1\bar{z'_1}+\cdots+z_n\bar{z'_n}$ is the Hermitean product on  $\mathbb{C}^n$. The Haar measure on $H^n$ coincides with the Lebesgue measure $dzdt$ on $\mathbb{C}^n\times\mathbb{R}$. It is easily to check that the inverse of the element $(z,s)$ is $(z,s)^{-1}=(-z,-s)$ and the unitary element is $(0,0)$.

The non-isotropic dilations on $H^n$ are defined by 
	\begin{equation*}
	\delta_a(z,s)=(az,a^2s),  \,  (a>0).
	\end{equation*}
The Haar measure on $H^n$ satisfies $d\delta_a(z,t) = a^Ndzdt$, where $N = 2n + 2$ is the homogeneous dimension of $H^n$. 

The Schwartz class of rapidly decreasing smooth function on $H^n$ is denoted by $\mathscr{S}(H^n)$ , and its dual, the space of tempered distribution on $H^n$, by $\mathscr{S}'(H^n)$. For $\Omega \subseteq H^n$, the space $C_c^\infty(\Omega)$ consists of smooth functions compactly supported on $\Omega$. We say a function $f$ is radial on $H^n$ if the value of $f(z,s)$ depends only on $|z|$ and $s$. 
We denote by $\mathscr{S}_{rad}(H^n)$ and $L_{rad}^p(H^n)$, $1\leqslant p \leq \infty$, the space of radial functions in $\mathscr{S}(H^n)$ and $L^p(H^n)$, respectively. In particular, the set of $L_{rad}^1(H^n)$ endowed with the convolution product 
\begin{equation*}
f*g(w)=\int_{H^n}f(wv^{-1})g(v)dv,\, w\in H^n
\end{equation*}
is a commutative algebra.

The left invariant vector fields on the Heisenberg group $H^n$ are defined by 
	\begin{equation*}
	\begin{aligned}
	X_j = \partial_{x_j} + 2y_j\partial_{s},  \\
	Y_j = \partial_{y_j} - 2x_j\partial_{s},
	\end{aligned}
	\end{equation*}
where $j=1,2,\cdots,n$. In terms of these vector fields we introduce the sub-Laplacian $\mathcal{L}$ by 
	\begin{equation*}
	\begin{aligned}
	\mathcal{L}&=-\underset{j=1}{\overset{n}{\sum}}(X_j^2 +Y_j^2) \\
	&=-\Delta_z-4|z|^2\partial_s^2+4\partial_s\sum_{j=1}^{n}(x_j\partial_{y_j}-y_j\partial_{x_j}).
	\end{aligned}
	\end{equation*}

\noindent\\[2mm]
\subsection{Spherical Fourier transform}
In this section, We discuss the spherical Fourier transform  on the Heisenberg group $H^n$.

The Laguerre polynomials of type $d\in\mathbb{N}$ and degree $m\in\mathbb{N}$ is defined by
	\begin{equation*}
	L_m^{(d)}(t)=\sum_{k=0}^{m}(-1)^k\left( \begin{array}{c} m+d\\m-k \end{array} \right)\frac{t^k}{k!},  \, t\geq0.
	\end{equation*}

\begin{definition} Let $f \in L_{rad}^1(H^n)$. For any $m \in \mathbb{N}$ and $\lambda \in \mathbb{R}$, we define the spherical Fourier transform of $f$ by
	\begin{equation*}
	\hat{f}(m,\lambda)=\left( \begin{array}{c} m+n-1\\m \end{array} \right)^{-1}
	\int_{H^n}e^{i\lambda s} f(z,s)L_m^{(n-1)}(2|\lambda||z|^2)e^{-|\lambda||z|^2}\,dzds.
	\end{equation*}	
\end{definition}
\noindent By a direct computation, we have  $\widehat{f*g}=\hat{f}\cdot\hat{g}$ for any $f,g\in L^1_{rad}(H^n)$.

The following inversion Fourier formula also holds.
\begin{proposition} For all $f \in \mathscr{S}_{rad}(H^n)$ such that
	\begin{equation*}
	\underset{m\in\mathbb{N}}{\sum}\left( \begin{array}{c} m+n-1\\m \end{array} \right) \int_{\mathbb{R}} |\hat{f}(\lambda,m)||\lambda|^d d\lambda <\infty,
	\end{equation*}
we have
	\begin{equation}\label{plancherel}	f(z,s)=\frac{2^{n-1}}{\pi^{n+1}}\underset{m\in\mathbb{N}}{\sum}\int_{\mathbb R} e^{-i\lambda\cdot s} \hat{f}(m,\lambda) L_m^{(n-1)}(2|\lambda||z|^2)e^{-|\lambda||z|^2}|\lambda|^n \,d\lambda,
	\end{equation}
where the sum is convergent in $L^{\infty}(H^n)$.
\end{proposition}
Moreover, if $f \in \mathscr{S}_{rad}(H^n)$, the functions $\mathcal{L} f$ is also in $\mathscr{S}_{rad}(H^n)$ and its spherical Fourier transform is given by
	\begin{equation*}
	\widehat{\mathcal{L} f}(m,\lambda)=b_m(\lambda)\hat{f}(m,\lambda),
	\end{equation*}
where $b_m(\lambda)=4(2m+n)|\lambda|$. The sub-Laplacian $\mathcal{L}$ is a positive self-adjoint operator densely defined on $L^2(H^n)$. So by the spectral theorem, for any bounded Borel function $h$ on $\mathbb{R}^+$, we have
	\begin{equation*}
	\widehat{h(\mathcal{L})f}(m,\lambda)=h(b_m(\lambda))\hat{f}(m,\lambda).
	\end{equation*}

\noindent\\[2mm]
\subsection{Homogeneous Besov spaces}
In this section, we shall recall the homogeneous Besov spaces related to the sub-Laplacian $\mathcal{L}$ on the Heisenberg group $H^n$ given in \cite{BGX2000}.
 
Let $R^* \in C_c^\infty(C_0)$ with $C_0=\{\tau \in \mathbb{R}:\frac{1}{2}\leq |\tau|\leq4\}$ be an even and positive function such that 
	\begin{equation*}
\underset{j\in \mathbb{Z}}{\sum}R^*(2^{-2j}\tau)=1, \quad \forall \tau \neq0.
	\end{equation*}
	
Since the restriction of $R^*$ on $\mathbb{R}^+$ is in $\mathscr{S}(\mathbb{R})$, by \cite{FMV2006},  we have the kernel $\varphi$ of the operator $R^*(\mathcal{L})\in\mathscr{S}_{rad}(H^n)$. In particular, 
	\begin{equation*}
	\widehat{\varphi}(m, \lambda)=R^*_m(\lambda):=R^*((2m+n)\lambda).
	\end{equation*}
By the inversion Fourier transform \eqref{plancherel}, we also have 
	\begin{equation*}
	\varphi(z,s)=\frac{2^{n-1}}{\pi^{n+1}}\sum_{m=0}^{\infty} \int_{\mathbb R} e^{-i\lambda\cdot s} R^*_m(\lambda) L_m^{(n-1)}(2|\lambda||z|^2)e^{-|\lambda||z|^2}|\lambda|^n \,d\lambda.
	\end{equation*}
	
For $j\in \mathbb{Z}$, we denote by $\varphi_j$ the kernel of the operator $R^*(2^{-2j}\mathcal{L})$. It can be easily calculated that
$\varphi_j(z,s)=2^{Nj}\varphi(\delta_{2^j}(z,s))$, $\forall (z,t)\in H^n$. The Littlewood-Paley projector on the Heisenberg group is define by
	\begin{equation*}
	\Delta_jf=f*\varphi_j,  \quad \forall f \in \mathscr{S}'(H^n).
	\end{equation*}
Obviously $\widehat{\varphi_j}(m,\lambda)=R^*_m(2^{-2j}|\lambda|):=R^j_m(\lambda)$. Let $\widetilde{\varphi_j}=\varphi_{j-1}+\varphi_j+\varphi_{j+1}$.  We have  $f*\varphi_j=f*\varphi_j*\widetilde{\varphi_j}$ holds true for all $f \in \mathscr{S}'(H^n)$.

By the spectral theorem, for any $f \in L^2(H^n)$, the following homogeneous Littlewood-Paley decompositon holds: 
	\begin{equation*}
	f=\sum_{j\in\mathbb{Z}}\Delta_j f \quad \text{in $L^2(H^n)$}.
	\end{equation*}
So
	\begin{equation}\label{infty}
	\|f\|_{L^\infty(H^n)}\leq\sum_{j\in\mathbb{Z}}\|\Delta_j f\|_{L^\infty(H^n)}, \, f\in L^2(H^n),
	\end{equation}
where both sides of \eqref{infty} are allowed to be infinite. 

In addition, it has been proved in \cite{FMV2006} that 
for any $\sigma\in\mathbb{R}$, $j\in\mathbb{Z}$, $1\leq p\leq \infty$ and $f\in \mathscr{S}'(H^n)$, then
	\begin{equation}\label{LPLP}
	\left\|\mathcal{L}^{\frac{\sigma}{2}}\Delta_jf\right\|_{L^p(H^n)} \leq C2^{j\sigma}\|\Delta_jf\|_{L^p(H^n)}.
	\end{equation}
\begin{definition} 
Let $1\leq p,r \leq \infty, \rho <N/p$. The homogeneous Besov space $\dot{B}^\rho_{p,r}(H^n)$ is defined as the set of distributions $f \in \mathscr{S}'(H^n)$ such that 
	\begin{equation*}
	\|f\|_{\dot{B}^\rho_{p,r}(H^n)}=\left(\underset{j\in \mathbb{Z}}{\sum}2^{j\rho r}\|\Delta_j\|_p^r\right)^{\frac{1}{r}}<\infty,
	\end{equation*}
and $f=\underset{j\in \mathbb{Z}}{\sum}\Delta_jf$ in $\mathscr{S}'(H^n)$.
\end{definition}
\begin{definition} 
Let $\rho <N/2$. The homogeneous Sobolev space $\dot{H}^\rho(H^n)$ is
	\begin{equation*}
	\dot{H}^\rho(H^n)=\dot{B}^0_{2,2}(H^n),
	\end{equation*}
which is equivalent to
	\begin{equation*}
	u\in\dot{H}^\rho(H^n)\Leftrightarrow\Delta^{\rho/2}u\in L^2(H^n).
	\end{equation*}
and the associated norms are of course equivalent.
\end{definition}
According to \cite{BGX2000}, we list some properties of the spaces $\dot{B}^\rho_{p,r}(H^n)$ in the following proposition.
\begin{proposition}\label{properties} Let $p,r\in [1,\infty]$ and $\rho<N/p$.
	
(i) The space $\dot{B}^\rho_{p,r}(H^n)$ is a Banach space with the norm $||\cdot||_{\dot{B}^\rho_{p,r}(H^n)}$;

(ii) the definition of $\dot{B}^\rho_{p,r}(H^n)$ does not depend on the choice of the function $R^*$ in the Littlewood-Paley decomposition;

(iii) for $-\frac{N}{p'}<\rho<\frac{N}{p}(H^n)$ the dual space of $\dot{B}^\rho_{p,r}(H^n)$ is $\dot{B}^{-\rho}_{p',r'}$;

(iv) for any $u\in \mathscr{S}'(H^n)$ and $\sigma>0$, then $u\in\dot{B}^\rho_{p,r}(H^n)$ if and only if $L^{\sigma/2}u\in\dot{B}^{\rho-\sigma}_{p,r}(H^n)$;

(v) for any $p_1,p_2\in[1,\infty]$, the continuous inclusion holds
\begin{equation*}
\dot{B}^{\rho_1}_{p_1,r}(H^n)\subseteq\dot{B}^{\rho_2}_{p_2,r}(H^n), \quad \frac{1}{p_1}-\frac{\rho_1}{N}=\frac{1}{p_2}-\frac{\rho_2}{N}, \quad \rho_1\geq\rho_2;
\end{equation*}

(vi) for all $p\in[2, \infty]$ we have the continuous inclusion $\dot{B}^0_{p,2}(H^n)\subseteq L^p(H^n)$;\\

(vii) $\dot{B}^0_{2,2}(H^n)=L^2(H^n)$.
\end{proposition}

\noindent\\[4mm]
\section{Technical Lemmas}
By the inversion Fourier formula \eqref{plancherel}, we may write $U_t\varphi_j$ explicitly into a sum of a list of oscillatory integrals. In order to estimate the oscillatory integrals, we recall the stationary phase lemma.
\begin{lemma}(see \cite{S1993})\label{phase} Let $g\in C^\infty([a,b])$ be real-valued such that
	\begin{equation*}
	|g''(x)|\geq \delta
	\end{equation*}
for any $x\in[a,b]$ with $\delta >0$. Then for any function $\psi \in C^\infty([a,b])$, there exists a constant $C$ which does not depend on $\delta, a, b, g$ or $\psi$, such that
	\begin{equation*}
	\left|\int_a^b e^{itg(x)}\psi(x)\,dx\right|\leq Ct^{-1/2}\left(\|\psi\|_\infty+\|\psi'\|_1\right).
	\end{equation*}
\end{lemma}

In order to prove the sharpness of the time decay dispersion of the solutions of some concrete equations in Section 5, we describe the asymptotic expansion of oscillating integrals.
\begin{lemma}(see \cite{S1993})\label{asymptotic} Suppose $k\geqslant2$, and $\phi(x_0)=\phi'(x_0)=\cdots=\phi^{k-1}(x_0)$, while $\phi^k(x_0)\not=0$. If $\psi$ is supported in a sufficiently neighborhood of $x_0$, then 
	\begin{equation*}
	I(t)=\left|\int_{\mathbb{R}}e^{it\phi(x)}\psi(x)\,dx\right| \sim |t|^{-1/k}\sum_{j=0}^{\infty}a_j|t|^{-j/k}, \quad \text{ as t}\rightarrow \infty.\\
	\end{equation*}
\end{lemma}

Besides, it will involve the Laguerre functions when we handle the oscillatory integrals. We need the following estimates.
\begin{lemma}(see \cite{H2005})\label{Laguerre} Consider the following:
	\begin{equation*}
	\left|(\tau \frac{d}{d\tau})^\alpha (L_m^{(d-1)}(\tau)e^{-\tau/2})\right|\leq C_{\alpha,d}(2m+d)^{d-1/4}
	\end{equation*}
for all $0\leq \alpha\leq d$.
\end{lemma}
\begin{remark}
In fact, for $0\leqslant \alpha \leqslant d-1$, we have a better estimate
	\begin{equation*}
	\left|(\tau \frac{d}{d\tau})^\alpha (L_m^{(d-1)}(\tau)e^{-\tau/2})\right|\leq C_{\alpha,d}(2m+d)^{d-1}.
	\end{equation*}
\end{remark}

Finally, we apply the following duality arguments to prove the Strichartz estimates for some concrete equations on the Heisenberg group. 
 
\begin{lemma}\label{equ} (see \cite{GV}) Let $H$ be a Hilbert space, $X$ a Banach space,$X^*$ the dual of $X$, and $D$ a vector space densely contained in $X$. Let $A\in \mathscr{L}_a(D,H)$ and let $A^*\in \mathscr{L}_a(H,D_a^*)$ be its adjoint, defined by
	\begin{equation*}
	\langle A^*v,f \rangle_D=\langle v,Af \rangle_H, \quad \forall f\in D, \quad \forall v \in H,
	\end{equation*}
where $\mathscr{L}_a(Y,Z)$ is the space of linear maps from a vector space $Y$ to a vector space $Z$, $D^*_a$ is the algebraic dual of $D$, $\langle \varphi,f\rangle_D$ is the pairing between $D^*_a$ and $D$, and $\langle,\rangle_H$ is the scalar product in $H$. Then the following three condition are equivalent:

(1) There exists $a$, $0\leqslant a \leqslant \infty$ such that for all $f \in D$,
	\begin{equation*}
	\|Af\| \leqslant a\|f\|_X.
	\end{equation*}
	
(2) $\mathfrak{R}(A^*)\subset X^*$, and there exists $a$, $0\leqslant a \leqslant \infty$, such that for all $v \in H$,
	\begin{equation*}
	\|A^*v\|_{X^*} \leqslant a\|v\|.
	\end{equation*}
	
(3) $\mathfrak{R}(A^*A)\subset X^*$,and there exists $a$, $0\leqslant a \leqslant \infty$, such that for all $f \in D$,
	\begin{equation*}
	\|A^*Af\|_{X^*} \leqslant a^2\|f\|_X,
	\end{equation*}
where $||\cdot||$ denote the norm in $H$. The constant $a$ is the same in all three parts. If one of those conditions is satisfied, the operators $A$ and $A^*A$ extend by continuity to bounded operators from $X$ to $H$ and from $X$ to $X^*$,respectively.
\end{lemma}
\begin{lemma}(see \cite{GV})\label{DuilatyXX} Let $H$,$D$ and two triplets$(X_i,A_i,a_i)$,$i=1,2$, satisfy any of the conditions of Lemma\ref{equ}. Then for all choices of $i,j=1,2$, $\mathfrak{R}(A_i^*A_j) \subset X_i^*$,and for all $f \in D$,
	\begin{equation*}
	\|A_i^* A_jf\|_{X_i^*} \leqslant a_i a_j\|f\|_{X_j}.
	\end{equation*}
\end{lemma}
\begin{lemma}\label{L1LW} (see \cite{GV}) Let $H$ be a Hilbert space, let $I$ be an interval of $\mathbb{R}$, let $X \subset S'(I \times \mathbb{R})$ be a Banach space, let $X$ be stable under time restriction, let $X$ and $A$ satisfy(any of) the conditions of Lemma \ref{equ}. Then the operator $A^*A$ is a bounded operator from $L_t^1(I,H)$ to $X^*$ and from $X$ to $L_t^\infty(I,H)$.
\end{lemma}

\noindent\\[4mm]
\section{Decay Estimates}
In this section, we will prove Theorem\ref{ResultTime}.

\noindent{\bf Proof of Theorem \ref{ResultTime}: }~ First, we prove (1). According to $\varphi_j$ introduced in section 2, we have
	\begin{equation*}
	\begin{aligned}
	U_t\varphi_j(z,s)&=\frac{2^{n-1}}{\pi^{n+1}}\sum_{m=0}^{\infty} \int_{-\infty}^{+\infty} e^{-i\lambda s} e^{it\phi(4(2m+n)|\lambda|)} \\ & \qquad \qquad \qquad \times R^j_m(\lambda) L_m^{(n-1)}(2|\lambda||z|^2)e^{-|\lambda||z|^2}|\lambda|^n \,d\lambda.
	\end{aligned}
\end{equation*}
By the transform of variable $\lambda'=2^{-2j}\lambda$, we get
	\begin{equation*}
	\begin{aligned}
	U_t\varphi_j(z,s)&=\frac{2^{n-1}}{\pi^{n+1}}\cdot2^{Nj}\sum_{m=0} ^\infty\int_{-\infty}^{+\infty} e^{-i2^{2j}\lambda' s} e^{it\phi((2m+n)2^{2j+2}|\lambda'|)} \\ 
	&\qquad \qquad \times R^*_m(\lambda') L_m^{(n-1)}( 2^{2j+1}|\lambda'||z|^2)e^{-2^{2j}|\lambda'||z|^2}|\lambda'|^n \,d\lambda' \\
	&=\frac{2^{n-1}}{\pi^{n+1}}\cdot2^{Nj}\sum_{m=0} ^\infty I_m(t,z,s).
	\end{aligned}
	\end{equation*}
Denote
	\begin{equation*}
	\begin{aligned}
	I_m(t,z,s)&=\int_{-\infty}^{+\infty} e^{-i2^{2j}\lambda' s} e^{it\phi((2m+n)2^{2j}+2|\lambda'|)}  \\
	& \times R^*_m(\lambda')L_m^{(n-1)}(2^{2j+1}|\lambda'||z|^2)e^{-2^{2j}|\lambda'||z|^2}|\lambda'|^n \,d\lambda'.
	\end{aligned}
	\end{equation*}
Let $s=ts'$, and $x=(2m+n)\lambda'$. We have
	\begin{equation*}
	\begin{aligned}
	I_m(t,z,s)&=\int_{C_0} e^{it\left(\phi( 2^{2j+2}|x|)-\frac{2^{2j}xs'}{2m+n}\right)}R^*(x) \\
	&\times L_m^{(n-1)}(\frac{ 2^{2j+1}|x||z|^2}{2m+n})e^{-\frac{2^{2j}|x||z|^2}{2m+n}}\frac{|x|^n}{(2m+n)^{n+1}} \,dx.
	\end{aligned}
	\end{equation*}	
It follows from Lemma \ref{Laguerre} that
	\begin{equation*}
	|I_m(t,z,s)| \leq C_n(2m+n)^{-n-1}\cdot (2m+n)^{n-1}\int_{C_0} R^*(x)|x|^n \,dx\leq C_n(2m+n)^{-2}.
	\end{equation*}
So
	\begin{equation*}
	|U_t\varphi_j(z,s)| \leq C_n2^{Nj} \sum_{m=0} ^\infty (2m+n)^{-2} \leq C_n2^{Nj}.
	\end{equation*}	
Finally we get
	\begin{equation}\label{res3-6}
	\|U_t \varphi_j\|_{L^\infty(H^n)}\leq C_n2^{Nj},
	\end{equation}	
which is the result of \eqref{res3-2}, as desired.	

Now we assume (H3) also holds. Let $g(x)=\phi(2^{2j+2}|x|)-\dfrac{2^{2j}xs'}{2m+n}$. By a direct computation, we obtain
	\begin{equation*}
	g''(x)= 2^{4j+4} \phi''(2^{2j+2}|x|).
	\end{equation*}
then $|g''(x)| \gtrsim 2^{2j\alpha_1}$ on the support of $R^*$. According to Lemma \ref{phase}, we have
	\begin{equation*}
	\begin{aligned}
	|I_m(t,z,s)| &\leq |t\cdot 2^{2j\alpha_1}|^{-\frac{1}{2}}\Biggl(\left\|R^*(x)L_m^{(n-1)}(\frac{ 2^{2j+1}|x||z|^2}{2m+n})e^{-\frac{2^{2j}|x||z|^2}{2m+n}}\frac{|x|^n}{(2m+n)^{n+1}}\right\|_{L^\infty(C_0)}  \\
	&+\int_{C_0}\left|\frac{d}{dx}\big(R^*(x)L_m^{(n-1)}(\frac{ 2^{2j+1}|x||z|^2}{2m+n})e^{-\frac{2^{2j}|x||z|^2}{2m+n}}\frac{|x|^n}{(2m+n)^{n+1}}\big)\right|\,dx\Biggr).
	\end{aligned}
	\end{equation*}	
Since $R^*\in C_c^\infty(C_0)$, by Lemma \ref{Laguerre}, it yields that
	\begin{equation*}
	\left\|R^*(x)L_m^{(n-1)}(\frac{ 2^{2j+1}|x||z|^2}{2m+n})e^{-\frac{2^{2j}|x||z|^2}{2m+n}}\frac{|x|^n}{(2m+n)^{n+1}}\right\|_{L^\infty(C_0)}  \leq C_n(2m+n)^{-2},
	\end{equation*}		
	\begin{equation*}
	\left|(\frac{d}{dx}R^*(x))L_m^{(n-1)}(\frac{ 2^{2j+1}|x||z|^2}{2m+n})e^{-\frac{2^{2j}|x||z|^2}{2m+n}}\frac{|x|^n}{(2m+n)^{n+1}}\right|\leq C_n(2m+n)^{-2},
	\end{equation*}
	\begin{equation*}
	\left|R^*(x)L_m^{(n-1)}(\frac{ 2^{2j+1}|x||z|^2}{2m+n})e^{-\frac{2^{2j}|x||z|^2}{2m+n}}\frac{1}{(2m+n)^{n+1}}\frac{d}{dx}(|x|^n)\right|\leq C_n(2m+n)^{-2},
	\end{equation*}
	\begin{equation*}
	\begin{aligned}
	&
	\left|R^*(x)\frac{d}{dx}(L_m^{(n-1)}(\frac{ 2^{2j+1}|x||z|^2}{2m+n})e^{-\frac{2^{2j}|x||z|^2}{2m+n}})\frac{|x|^n}{(2m+n)^{n+1}}\right|  \\
	=& R^*(x)\cdot\left|y\frac{d}{dy}(L_m^{(n-1)}(2y)e^{-y})\right|_{y=\frac{ 2^{2j}|x||z|^2}{2m+n}}\cdot|x|^{n-1}\cdot \frac{1}{(2m+n)^{n+1}} \\
	\leq & C_n(2m+n)^{-2}.
	\end{aligned}
	\end{equation*}
So we get 
	\begin{equation*}
		|I_m(t,z,s)| \leq C_n|t|^{-\frac{1}{2}} \cdot 2^{-j\alpha_1}(2m+n)^{-2},
	\end{equation*}
and 
	\begin{equation*}
	|U_t\varphi_j(z,s)| \leq C_n|t|^{-\frac{1}{2}}\cdot 2^{(N-\alpha_1)j}\sum_{m=0} ^\infty (2m+n)^{-2}  \leq C_n|t|^{-\frac{1}{2}}\cdot 2^{(N-\alpha_1)j}.
	\end{equation*}	
Then
	\begin{equation}\label{res3-7}
	\|U_t \varphi_j\|_{L^\infty(H^n)} \leq C_n|t|^{-\frac{1}{2}}\cdot 2^{(N-\alpha_1)j}.
	\end{equation}	
By an interpolation between \eqref{res3-6} and \eqref{res3-7}, we obtain 
	\begin{equation*}
	\|U_t \varphi_j\|_\infty \leq C_n|t|^{-\frac{\theta}{2}}\cdot 2^{j(N-\alpha_1\theta)},\,\text{for }0\leq \theta \leq 1,
	\end{equation*}	
which completes the proof of (1).

 The proof of (2) is similar to (1) and we omit the details.\\
 	\qed	
	

\noindent\\[4mm]
\section{Applications}
In this section, we prove the Strichartz estimates for some concrete equations on the Heisenberg group by using Theorem \ref{ResultTime}. 

\subsection{The Fractional Schr\"{o}dinger Equation}

\indent First, we consider the fractional Schr\"{o}dinger equation ($0<\alpha<1$)
    \begin{equation}\label{FSchrEqu}
	\begin{cases}
	i\partial_tu+\mathcal{L}^\alpha u=f,\\
	u|_{t=0}=u_0.
	\end{cases}
	\end{equation}
The solution is given by
	\begin{equation*}\label{solution1}
	u(t)=S_tu_0-i\int_0^tS_{t-\tau}f(\tau)\,d\tau,
	\end{equation*}
where $S_t=e^{it\mathcal{L}^\alpha}$ and it corresponds to the case when $\phi(r)=r^\alpha$. By a simple calculation, 
	\begin{equation*}
	\phi'(\lambda)=\alpha r^{\alpha-1}, \quad
	\phi''(\lambda)=\alpha(\alpha-1)r^{\alpha-2}.
	\end{equation*}
We see that $\phi$ satisfies (H1)-(H4) with $m_1=\alpha_1=m_2=\alpha_2=\alpha$. So according to Theorem\ref{ResultTime}, we know 
	\begin{equation*}
	\|S_t \varphi_j\|_{L^\infty (H^n)} \leq C_n|t|^{-\frac{\theta}{2}}2^{j(N-\alpha\theta)}, \quad \forall j \in \mathbb{Z}.
	\end{equation*}
We notice that when $\theta=1$ the fractional Schr\"{o}dinger semigroup has the best decay in time
	\begin{equation}\label{SchTime1}
	\|S_t \varphi_j\|_{L^\infty (H^n)} \leq C_n|t|^{-\frac{1}{2}}2^{j(N-\alpha)}, \quad \forall j \in \mathbb{Z}.
	\end{equation}

\begin{theorem}
If $u$ is the solution of the free Schr\"{o}dinger equation \eqref{FSchrEqu} (i.e. with $f=0$). For $0<\alpha<1$, there exists a constant $C>0$ which does not depend on $u_0$ or $t$ such that
	\begin{equation*}
	\|u(t)\|_{L^\infty (H^n)}\leq C|t|^{-1/2}\|u_0\|_{\dot{B}^{N-\alpha}_{1,1}(H^n)},
	\end{equation*}
and the time decay is sharp in time.
\end{theorem}
{\bf Proof.}~According to the property of the Fourier transform, for any $f\in \mathscr{S}(H^n)$, we have
	\begin{equation*}
	\begin{aligned}
	\Delta_j(S_tf)&=S_t(\Delta_jf) \\
	&=S_t(f*\varphi_j)  \\
	&=S_t(f*\varphi_j*\widetilde{\varphi_j})  \\
	&=f*\varphi_j*(S_t\widetilde{\varphi_j}).
	\end{aligned}
	\end{equation*} 
By the Young inequality and \eqref{SchTime1}, we get
	\begin{equation}\label{L1LWinequality}
	\begin{aligned}
	\left\|\Delta_j(S_tf)\right\|_{L^\infty (H^n)} &\leq \left\|f*\varphi_j\right\|_{L^1(H^n)}\cdot \left\|S_t\widetilde{\varphi_j}\right\|_{L^\infty (H^n)}   \\ 
	&\leq C_n|t|^{-\frac{1}{2}}2^{j(N-\alpha)}\|\Delta_jf\|_{L^1(H^n)}.
	\end{aligned}
	\end{equation}
Finally, we have the following results:
	\begin{equation*}
	\begin{aligned}
	\|u(t)\|_{L^\infty (H^n)}&=\|S_tu_0\|_{L^\infty (H^n)}\\
	&=\|\sum_{j\in \mathbb{Z}}\Delta_j(S_tu_0)\|_{L^\infty (H^n)} \\
	&\leq \sum_{j\in \mathbb{Z}}\left\|\Delta_j(S_tu_0)\right\|_{L^\infty (H^n)}\\
	&\leq \sum_{j\in \mathbb{Z}}C|t|^{-\frac{1}{2}}2^{j(N-\alpha)}\|\Delta_j u_0\|_{L^1(H^n)} \\
	&=C|t|^{-\frac{1}{2}}\|u_0\|_{\dot{B}^{N-\alpha}_{1,1}(H^n)}.
	\end{aligned}
	\end{equation*}

Now let us give an example to explain the sharpness of the decay in time $t$. Let $Q\in C_0^\infty(D_0)$ with $Q(1)=1$, where $D_0$ is a small neighborhood of 1 such that $0\notin D_0$. Then
	\begin{equation*}
	\begin{aligned}
	\widehat{u_0}(m,\lambda)&=Q(|\lambda|)\delta_{m0},   \\
	u_0(z,s)&=\int_{D_0}e^{-i\lambda s}Q(\lambda)\lambda^n \,d\lambda.
	\end{aligned}
	\end{equation*}
Consider the free fractional Schr\"{o}dinger equation (\ref{FSchrEqu}). We have
	\begin{equation*}	
	u(z,s,t)=e^{it\mathcal{L}^\alpha}u_0=C\int_{D_0}e^{-i\lambda s}e^{it(4\lambda n)^\alpha }Q(\lambda)e^{-\lambda|z|^2}\lambda^n \,d\lambda.
	\end{equation*}
In particular,
	\begin{equation*}
	u(0,ts,t)=C\int_{D_0}e^{it(-\lambda s+4^\alpha\lambda^\alpha n^\alpha)}Q(\lambda)\lambda^n \,d\lambda.
	\end{equation*}
This oscillating integral has a phase $\Psi(\lambda)=-\lambda s+4^\alpha\lambda^\alpha n^\alpha$ with a unique critical point $\lambda_0=(\dfrac{s}{\alpha4^\alpha n^\alpha})^{\frac{1}{\alpha-1}}$, which is not degenerate. By Lemma \ref{asymptotic}, for $s_0=\alpha4^\alpha n^\alpha$, we get
	\begin{equation*}
	u(0,ts_0,t) \sim C|t|^{-1/2}.
	\end{equation*}
\qed

\begin{proposition}\label{Intermediate} Let $0<\alpha<1$. For $i=1,2$, let $p_i, r_i\in[2,\infty]$ and $\rho_i\in\mathbb{R}$ such that
	\begin{equation*}
	\begin{aligned}
	&(i)\,\dfrac{2}{p_i}+\dfrac{1}{r_i} \leq \dfrac{1}{2},   \\
	&(ii)\,\frac{N}{r_i}+\frac{2\alpha}{p_i}-\rho_i=\frac{N}{2},
	\end{aligned}
	\end{equation*}
Let $\bar{p}_i, \bar{r}_i$ denote the conjugate exponent of $p_i, r_i$ for $i=1,2$. Then the following estimates are satisfied:
	\begin{equation}\label{SEresult1}
	\|S_tu_0\|_{L^{p_1}(\mathbb{R}, \dot{B}_{r_1,2}^{\rho_1}(H^n))}\leq C\|u_0\|_{L^2(H^n)},
	\end{equation}
	\begin{equation}\label{SEresult2}
	\left\|\int_0^tS_{t-\tau}f(\tau)\,d\tau\right\|_{L^{p_1}[0,T], \dot{B}^{\rho_1}_{r_1,2}(H^n))}\leq C\|f\|_{L^{\bar{p}_2}([0,T],\dot{B}_{\bar{r}_2,2}^{-\rho_2}(H^n))},
	\end{equation}
where the constant $C>0$ does not depend on $u_0$, $f$ or $T$.
\end{proposition}

{\bf Proof.} Firstly, according to Plancherel theorem, it is easy to obtain
	\begin{equation*}\label{plancherel1}
	\|\Delta_j(S_tf)\|_{L^2(H^n)}=\|S_t(\Delta_jf)\|_{L^2(H^n)}=\|\Delta_jf\|_{L^2(H^n)}.
	\end{equation*}
By an interpolation between the above equation and \eqref{L1LWinequality}, for any $2\leqslant r \leqslant \infty$, we have
	\begin{equation*}\label{LrLr1}
	\|\Delta_j(S_tf)\|_{L^r(H^n)} \leq C|t|^{-(1/2-1/r)}2^{2j(N-\alpha)(1/2-1/r)}\|\Delta_jf\|_{L^{\bar{r}}(H^n)}.
	\end{equation*}
Let $\gamma(r)=1/2-1/r$, $\beta(r)=(N-\alpha)(1/2-1/r)$. We have
	\begin{equation*}\label{LrLr}
	\|\Delta_j(S_tf)\|_{L^r(H^n)} \leq C|t|^{-\gamma(r)}2^{2j\beta(r)}\|\Delta_jf\|_{L^{\bar{r}}(H^n)}.
	\end{equation*}
Then 
	\begin{align}\label{BesovEstimate}
	\|S_tf\|_{\dot{B}_{r,2}^{-\beta(r)}}&=\left(\sum_{j=0}2^{-2j\beta(r)}\left\|\Delta_j(S_tf)\right\|^2_{L^r(H^n)}\right)^{1/2}    \nonumber   \\
	&\leq C\left(\sum_{j=0}2^{2j\beta(r)}|t|^{-2\gamma(r)}\left\|\Delta_jf\right\|^2_{L^{\bar{r}}(H^n)}\right)^{1/2}         \nonumber    \\
	& \leq C|t|^{-\gamma(r)}\left\|f\right\|_{\dot{B}_{\bar{r},2}^{\beta(r)}}.
	\end{align}
To prove \eqref{SEresult1}, it is sufficient to prove it for any $g\in L^2(H^n)$, $\psi \in L^{\bar{p}_1}(\mathbb{R},\dot{B}_{\bar{r}_1,2}^{-\rho_1}(H^n))$, the following holds:
	\begin{equation*}\label{L2LP}
	\left|\langle S_tg,\psi\rangle_{L^2(\mathbb{R},L^2(H^n))}\right| \leq C\|g\|_{L^2(H^n)}\|\psi\|_{L^{\bar{p}_1}(\mathbb{R},\dot{B}^{-\rho_1}_{\bar{r}_1,2})}.
	\end{equation*}
By density, we only need to prove for any $g\in L^2(H^n)$, \eqref{L2LP} holds for any $\psi \in \mathscr{S}(\mathbb{R},L^2(H^n))$.  

First, we have 
	\begin{equation}\label{L2Es}
	\begin{aligned}
	\int_{\mathbb{R} \times H^n}S_t g(z,s)\bar{\psi}(t,&z,s)\,dtdzds=\langle g,\int_{\mathbb{R}}S_t\psi
	(t)\,dt\rangle_{L^2(H^n)}\\
	&\leq \|g\|_{L^2(H^n)}\left\|\int_{\mathbb{R}}S_t\psi(t)\,dt\right\|_{L^2(H^n)} \\
	&=\|g\|_{L^2(H^n)}\left(
	\int_{\mathbb{R}}\int_{\mathbb{R}}\langle S_{t-\tau}\psi(t),\psi(\tau)\rangle_{L^2(H^n)} \,dtd\tau\right)^{1/2}.
	\end{aligned}
	\end{equation}
We can easily obtain
	\begin{equation*}\label{BesovTi1}
	\langle S_{t-\tau}\psi(t),\psi(\tau)\rangle_{L^2(H^n)} \leq \|S_{t-\tau}\psi(t)\|_{\dot{B}_{r,2}^{-\beta(r)}(H^n)}\|\psi(\tau)\|_{\dot{B}_{\bar{r},2}^{\beta(r)}(H^n)}.
	\end{equation*}
By \eqref{BesovEstimate} it yields that
	\begin{equation*}\label{BesovTi2}
	\left\|S_{t-\tau}\psi
	(t)\right\|_{\dot{B}_{r,2}^{-\beta(r)}(H^n)} \leq C|t-\tau|^{-\gamma(r)}\|\psi(t)\|_{\dot{B}_{\bar{r},2}^{\beta(r)}(H^n)}.
	\end{equation*}
Thus it holds
	\begin{align*}
	&\left(\int_{\mathbb{R}}\int_{\mathbb{R}}\langle S_{\tau-t}\psi(t),\psi(\tau)\rangle_{L^2(H^n)} \,dtd\tau\right)^{1/2} \\
	&\leq \left(\int_{\mathbb{R}}\int_{\mathbb{R}}\frac{\|\psi(t)\|_{\dot{B}_{\bar{r},2}^{\beta(r)}(H^n)}\|\psi(\tau)\|_{\dot{B}_{\bar{r},2}^{\beta(r)}(H^n)}}{|t-\tau|^{\gamma(r)}}\,dtd\tau\right)^{1/2}.
	\end{align*}	
By the Hardy-Littlewood-Sobolev inequality, we have
	\begin{equation*}\label{HLS}
	\begin{aligned}
	&\left(\int_{\mathbb{R}}\int_{\mathbb{R}}\frac{\|\psi(t)\|_{\dot{B}_{\bar{r},2}^{\beta(r)}(H^n)}\|\psi(\tau)\|_{\dot{B}_{\bar{r},2}^{\beta(r)}(H^n)}}{|t-\tau|^{\gamma(r)}}\,dtd\tau\right)^{1/2}  \\
	&=\left(\int_{\mathbb{R}}\int_{\mathbb{R}}\frac{\|\psi(t)\|_{\dot{B}_{\bar{r},2}^{\beta(r)}(H^n)}\|\psi(\tau)\|_{\dot{B}_{\bar{r},2}^{\beta(r)}(H^n)}}{|t-\tau|^{1-(1-\gamma(r))}}\,dtd\tau\right)^{1/2}  \\
	&\leq C\|\psi(t)\|_{L^{\bar{p}_1}(\mathbb{R},\dot{B}^{\beta(r)}_{\bar{r},2}(H^n))}.
	\end{aligned}
	\end{equation*}
where $\dfrac{2}{p_1}+\dfrac{1}{r}=\dfrac{1}{2}$, $\rho=-\beta(r)$. It follows from the above equation and \eqref{L2Es} that
	\begin{equation*}
	\left|\langle S_tg,\psi\rangle_{L^2(\mathbb{R},L^2(H^n))}\right| \leq C\|g\|_{L^2(H^n)}\|\psi\|_{L^{\bar{p}_1}(\mathbb{R},\dot{B}^{\beta(r)}_{\bar{r}_1,2}(H^n))}.
	\end{equation*}
Thus \eqref{SEresult1} holds.  

So far, we have proved \eqref{SEresult1} for the extreme case when $\dfrac{2}{p_1}+\dfrac{1}{r}=\dfrac{1}{2}$, $\dfrac{N}{r}+\dfrac{2\alpha}{p_1}-\rho=\dfrac{N}{2}$. For $\dfrac{2}{p_1}+\dfrac{1}{r_1}\leq\dfrac{1}{2}$ and $\dfrac{N}{r_1}+\dfrac{2\alpha}{p_1}-\rho_1=\dfrac{N}{2}$, we have $r_1\geq r$ and $\rho_1\leq\rho$. By Proposition \ref{properties}, we have $\dot{B}^\rho_{r,2}(H^n) \subset \dot{B}^{\rho_1}_{r_1,2}(H^n)$ and \eqref{SEresult1} also holds true.

According to \eqref{SEresult1}, Lemma \ref{equ} and Lemma \ref{DuilatyXX} we can immediately get
	\begin{equation*}
	\left\|\int_0^tS_{t-\tau}f(\tau)\,d\tau\right\|_{L^{p_1}([0,T], \dot{B}^{\rho_1}_{r_1,2}(H^n))}\leq C\|f\|_{L^{\bar{p}_2}([0,T],\dot{B}_{\bar{r}_2,2}^{-\rho_2}(H^n))}.
	\end{equation*}
Thus \eqref{SEresult2} holds.
In addition, by Lemma \ref{L1LW}, we also have
	\begin{equation}\label{LP1L1}
	\left\|\int_{0}^{t}S_{t-\tau}f(\tau)d\tau\right\|_{L^{p_1}(I,\dot{B}^{\rho_1}_{r_1,2}(H^n))} \leq C\|f\|_{L^1(I,L^2(H^n))}.
	\end{equation}
	\begin{equation*}
	\left\|\int_{0}^{t}S_{t-\tau}f(\tau)d\tau\right\|_{L^\infty(I,L^2(H^n))} \leq C\left\|\int_{0}^{t}S_{t-\tau}f(\tau)d\tau\right\|_{L^{\bar{p}_1}(I,\dot{B}^{-\rho_1}_{\bar{r}_1,2}(H^n))}.
	\end{equation*}\\
	\qed
	
\begin{theorem}
Under the same hypotheses as in Theorem \ref{Intermediate}, the solution $u$ of
the fractional Schr\"{o}dinger equation \eqref{FSchrEqu} satisfies the following estimate
	\begin{equation*}
	\|u\|_{L^{p_1}([0,T], \dot{B}^{\rho_1}_{r_1,2}(H^n))} \leq C\left(\|u_0\|_{L^2(H^n)}+\|f\|_{L^{\bar{p}_2}([0,T],\dot{B}_{\bar{r}_2,2}^{-\rho_2}(H^n))}\right),
	\end{equation*}
where the constant $C>0$ does not depend on $u_0$, $f$ or $T$.
\end{theorem}

Next, by the embedding relation between the homogeneous Besov space and Lebesgue space, we prove the Strichartz inequalities on Lebesgue spaces.
\begin{corollary}\label{SLebesgue} Let $0<\alpha<1$ and $u$ be the solution of the fractional Schr\"{o}dinger equation \eqref{FSchrEqu}. Suppose $p\in[2N-2\alpha,+\infty)$ and $q$ satisfy
	\begin{equation*}
	\dfrac{2\alpha}{p}+\dfrac{N}{q}=\dfrac{N}{2}-1.
	\end{equation*}
Then we have the following estimate
	\begin{equation*}
	\|u\|_{L^p([0,T],L^q(H^n))} \leq C\left(\|u_0\|_{\dot{H}^1(H^n)}+\|f\|_{L^1([0,T],\dot{H}^1(H^n))}\right),
	\end{equation*}
	where the constant $C>0$ does not depend on $u_0$, $f$ or $T$.
\end{corollary}
{\bf Proof.} Under the same hypotheses as in Proposition \ref{Intermediate} and according to \eqref{SEresult1} and \eqref{LP1L1}, we have
	\begin{equation*}	
	\begin{aligned}
	&\quad\|u\|_{L^p([0,T],\dot{B}^{\rho+1}_{r,2}(H^n))} \\&\leq \left\|S_tu_0\right\|_{L^p([0,T],\dot{B}^{\rho+1}_{r,2}(H^n))} +\left\|\int_{0}^{t}S_{t-s}f(s)ds\right\|_{L^p([0,T],\dot{B}^{\rho+1}_{r,2}(H^n))}.  \\
	 &\leq \|u_0\|_{\dot{H}^1(H^n)}+\|f\|_{L^1([0,T],\dot{H}^1(H^n))}.
	\end{aligned}
	\end{equation*}
Combining with the fact $\dot{B}^{\rho+1}_{p,r}(H^n)\subset L^q(H^n)$ where $0\leq\rho+1<N/r$ and $\dfrac{1}{q}=\dfrac{1}{r}-\dfrac{\rho+1}{N}$, we have 
	\begin{equation*}
	\|u\|_{L^p([0,T],L^q(H^n))} \leq C\left(\|u_0\|_{\dot{H}^1(H^n)}+\|f\|_{L^1([0,T],\dot{H}^1(H^n))}\right).
	\end{equation*}
	\qed
\begin{remark}
	It has been proved in \cite{BGX2000} that the Schr\"{o}dinger equation involving the sub-Laplacian on the Heisenberg group is not dispersive and for this reason the classical method to obtain Strichartz inequalities from dispersive estimates does not work. Our approach inspired by \cite{GPW2008} enables us to obtain dispersive estimates and also Strichartz inequalites for fractional Schr\"{o}dinger equations.
\end{remark}

\noindent\\[2mm]
\subsection{The Fractional Wave Equation}
We now consider the fractional wave equation ($0<\alpha<2$)
	\begin{equation}\label{FWaveEqu}
	\begin{cases}
	\partial_t^2u+\mathcal{L}^\alpha u=f \in L^1([0,T],L^2(H^n)),\\
	u|_{t=0}=u_0 \in \dot{H}^\alpha(H^n),  \\
	\partial_tu|_{t=0}=u_1 \in L^2(H^n),
	\end{cases}
	\end{equation}
The general solution of \eqref{FWaveEqu} has the form 
	\begin{equation*}
	u=v+w.
	\end{equation*}
where $v$ is a solution of \eqref{FWaveEqu} with $f=0$ and $w$ is the solution of \eqref{FWaveEqu} with $u_0=u_1=0$. They are classically given by
	\begin{equation}\label{Vsolution}
	v(t)=\dfrac{dA_t}{dt}u_0+A_tu_1,
	\end{equation}
	\begin{equation}\label{Wsolution}
	w(t)=\int_{0}^{t}A_{t-\tau}f(\tau)d\tau,
	\end{equation}
where 
	\begin{equation*}
	A_t=\dfrac{sin(t\mathcal{L}^{\frac{\alpha}{2}})}{\mathcal{L}^{\frac{\alpha}{2}}},
	\quad  \dfrac{dA_t}{dt}=cos(t\mathcal{L}^{\frac{\alpha}{2}}).
	\end{equation*}
So we naturally introduce the operator $W_t=e^{it\mathcal{L}^{\frac{\alpha}{2}}}$, which corresponds to the case when $\phi(r)=r^{\frac{\alpha}{2}}$. By a
simple calculation,
	\begin{equation*}
	\phi'(r)=\frac{\alpha}{2}\cdot r^{\frac{\alpha}{2}}, \quad
	\phi''(\lambda)
	=\frac{\alpha}{2}(\frac{\alpha}{2}-1)r^{\frac{\alpha}{2}-2}.
	\end{equation*}
It yields that $\phi$ satisfies (H1)-(H4) with $m_1=\alpha_1=m_2=\alpha_2=\dfrac{\alpha}{2}$. According to Theorem\ref{ResultTime}, we have
	\begin{equation}\label{WTimeE}
	\|W_t \varphi_j\|_{L^\infty (H^n)} \leq C_n|t|^{-\frac{1}{2}}2^{j(N-\frac{\alpha}{2})}, \quad \forall j\in \mathbb{Z}.
	\end{equation}
\begin{theorem}
Let $u_0\in\dot{B}^{N-\frac{\alpha}{2}}_{1,1}(H^n)$ and $u_1\in\dot{B}^{N-\frac{3}{2}\alpha}_{1,1}(H^n)$. Suppose $v$ is the solution of the free fractional wave equation \eqref{FWaveEqu} (i.e. with $f=0$). Then there exists a constant $C>0$, which does not depend on $u_0$ or $u_1$ such that
	\begin{equation}\label{WDecay}
	\|v(t)\|_{L^\infty (H^n)}\leq C|t|^{-1/2}\left(\|u_0\|_{\dot{B}^{N-\frac{\alpha}{2}}_{1,1}(H^n)}+\|u_1\|_{\dot{B}^{N-\frac{3}{2}\alpha}_{1,1}(H^n)}\right),
	\end{equation}
and the time decay is sharp in time.
\end{theorem}
{\bf Proof.}~Because of \eqref{Vsolution}, we can write 
	\begin{equation*}
	v(t)=\frac{dA_t}{dt}u_0+A_tu_1  
	=\frac{W_t+W_{-t}}{2}u_0+\mathcal{L}^{-\frac{\alpha}{2}}\frac{W_t-W_{-t}}{2i}u_1.	
	\end{equation*}
Therefore, it suffices to prove that $\|W_tu_0\|_{L^\infty (H^n)}$ and  $\|\mathcal{L}^{-\frac{\alpha}{2}}W_tu_1\|_{L^\infty (H^n)}$ can be controlled by the right side of \eqref{WDecay}. 

According to \eqref{WTimeE} and the Young inequality, we get
	\begin{equation*}
	\begin{aligned}
	\|W_tu_0\|_{L^\infty (H^n)}&=\|\sum_{j\in \mathbb{Z}}\Delta_j(W_tu_0
	)\|_{L^\infty (H^n)} \leq \sum_{j\in \mathbb{Z}}\|\Delta_j(W_tu_0
	)\|_{L^\infty (H^n)}  \\
	&=\sum_{j\in \mathbb{Z}}\|W_t(\Delta_ju_0)
	\|_{L^\infty (H^n)}=\sum_{j\in \mathbb{Z}}\|W_t(u_0*\varphi_j)
	\|_{L^\infty (H^n)} \\
	&=\sum_{j\in \mathbb{Z}}\|W_t(u_0*\varphi_j*\widetilde{\varphi_j})
	\|_{L^\infty (H^n)} =\sum_{j\in \mathbb{Z}}\|u_0*\varphi_j*W_t\widetilde{\varphi_j}
	\|_{L^\infty (H^n)} \\
	&\leq \sum_{j\in \mathbb{Z}}\|u_0*\varphi_j\|_{L^1 (H^n)}\|W_t\widetilde{\varphi_j}\|_{L^\infty (H^n)}\\
	&\leq \sum_{j\in \mathbb{Z}}C|t|^{-\frac{1}{2}}2^{j(N-\frac{\alpha}{2})}\|\Delta_ju_0\|_{L^1(H^n)}\\
	&=C|t|^{-\frac{1}{2}}\|u_0\|_{\dot{B}^{N-\frac{\alpha}{2}}_{1,1}(H^n)},
	\end{aligned}
	\end{equation*} 
and 
	\begin{equation*}
	\left\|\mathcal{L}^{-\frac{\alpha}{2}}W_tu_1\right\|_{L^\infty (H^n)}=\left\|W_t(\mathcal{L}^{-\frac{\alpha}{2}}u_1)\right\|_{L^\infty (H^n)} \leq C|t|^{-\frac{1}{2}}\left\|\mathcal{L}^{-\frac{\alpha}{2}}u_1\right\|_{\dot{B}^{N-\frac{\alpha}{2}}_{1,1}}. 
	\end{equation*}
It follows from \eqref{LPLP} that
	\begin{equation*}
	\left\|\Delta_j(\mathcal{L}^{-\frac{\alpha}{2}}u_1)\right\|_{L^1(H^n)}=\left\|(\mathcal{L}^{-\frac{\alpha}{2}}(\Delta_ju_1))\right\|_{L^1(H^n)} \leq
	C2^{-j\alpha}\|\Delta_ju_1\|_{L^1(H^n)}.
	\end{equation*}
Hence, we have 
	\begin{equation*}
	\left\|\mathcal{L}^{-\frac{\alpha}{2}}W_tu_1\right\|_{L^\infty (H^n)} \leq C|t|^{-\frac{1}{2}}\|u_1\|_{\dot{B}^{N-\frac{3}{2}\alpha}_{1,1}(H^n)}.
	\end{equation*}

Next, we also give an example to illustrate the sharpness of the decay in time $t$.  For the same $Q$, $D_0$ and $u_0$ as in the case of the fractional Sch\"{o}dinger equation, we consider the fractional wave equation (\ref{FWaveEqu}) when $u_1=0$, $f=0$. We have
	\begin{equation*}	
	u(z,s,t)=e^{it\mathcal{L}^\frac{\alpha}{2}}u_0=C\int_{D_0}e^{-i\lambda s}e^{it(4\lambda n)^\frac{\alpha}{2} }Q(\lambda)e^{-\lambda|z|^2}\lambda^n \,d\lambda.
	\end{equation*}
In particular,
	\begin{equation*}
	u(0,ts,t)=C\int_{D_0}e^{it(-\lambda s+(4\lambda n)^\frac{\alpha}{2})}Q(\lambda)\lambda^n \,d\lambda.
	\end{equation*}
This oscillating integral has a phase $\Psi(\lambda)=-\lambda s+(4\lambda n)^\frac{\alpha}{2}$ with a unique critical point $\lambda_0=\dfrac{1}{4n}\left(\dfrac{2\alpha n}{s}\right)^{\frac{2}{2-\alpha}}$, which is not degenerate. By Lemma \ref{asymptotic}, for $s_0=\dfrac{\alpha(4n)^\frac{\alpha}{2}}{2}$, it yields that
	\begin{equation*}
	u(0,ts_0,t) \sim C|t|^{-1/2}.
	\end{equation*}
\\
	\qed
	
As Propostion \ref{Intermediate}, we have the following Strichartz estimates.
\begin{proposition}\label{WOStrEst}
Suppose $0<\alpha<2$. For $i=1,2$, let $p_i, r_i\in[2,\infty]$ and $\rho_i\in\mathbb{R}$ such that
	\begin{equation*}
	\begin{aligned}
	&(i)\,\dfrac{2}{p_i}+\dfrac{1}{r_i} \leq \dfrac{1}{2},   \\
	&(ii)\,\frac{N}{r_i}+\frac{\alpha}{p_i}-\rho_i=\frac{N}{2},
	\end{aligned}
	\end{equation*}
Then the following estimates are satisfied:
	\begin{equation}\label{WEresult1}
	\|W_tg\|_{L^{p_1}(\mathbb{R}, \dot{B}_{r_1,2}^{\rho_1}(H^n))}\leq C\|g\|_{L^2(H^n)},
	\end{equation}
	\begin{equation}\label{WEresult2}
	\left\|\int_0^tW_{t-\tau}f(\tau)\,d\tau\right\|_{L^{p_1}([0,T], \dot{B}^{\rho_1}_{r_1,2}(H^n))}\leq C\|f\|_{L^{\bar{p}_2}([0,T],\dot{B}_{\bar{r}_2,2}^{-\rho_2}(H^n))},
	\end{equation}
where the constant $C>0$ does not depend on $u_0$, $f$ or $T$.
\end{proposition}

Now we describe the Strichartz inequalities for the solution of the fractional wave equation \eqref{FWaveEqu}.
\begin{theorem}
Suppose $0<\alpha<2$. For $i=1,2$, let $p_i, r_i\in[2,\infty]$ and $\rho_i\in\mathbb{R}$ such that
	\begin{equation*}
	\begin{aligned}
	&(i)\,\dfrac{2}{p_i}+\dfrac{1}{r_i} \leq \dfrac{1}{2},   \\
	&(ii)\,\frac{N}{r_1}+\frac{\alpha}{p_1}-\rho_1=\frac{N}{2}-\alpha,   \\
	&(iii)\,\frac{N}{r_2}+\frac{\alpha}{p_2}-\rho_2=\frac{N}{2},
	\end{aligned}
	\end{equation*}
Then the following estimates are satisfied:
	\begin{equation}\label{WResult1}
	\|v\|_{L^{p_1}(\mathbb{R},\dot{B}_{r_1,2}^{\rho_1}(H^n))}+\|\partial_tv\|_{L^{p_1}(\mathbb{R},\dot{B}_{r_1,2}^{\rho_1-\alpha}(H^n))} \leq C\left(\|u_0\|_{\dot{H}^\alpha(H^n)}+\|u_1\|_{L^2(H^n)}\right),
	\end{equation}
	\begin{equation}\label{WResult2}
	\|w\|_{L^{p_1}([0,T],\dot{B}_{r_1,2}^{\rho_1}(H^n))}+\|\partial_tw\|_{L^{p_1}([0,T],\dot{B}_{r_1,2}^{\rho_1-\alpha}(H^n))} \leq C\|f\|_{L^{\bar{p}_2}([0,T],\dot{B}_{\bar{r}_2,2}^{-\rho_2}(H^n))},
	\end{equation}
where the constant $C>0$ does not depend on $u_0$, $f$ or $T$.
\end{theorem}
{\bf Proof.}~First, we can write
	\begin{equation*}
	v=\frac{dA_t}{dt}u_0+A_tu_1  
	=\frac{W_t+W_{-t}}{2}u_0+\mathcal{L}^{-\frac{\alpha}{2}}\frac{W_t-W_{-t}}{2i}u_1,	
	\end{equation*}
and
	\begin{equation*}
	\partial_tv=A_t\mathcal{L}^\alpha u_0+\frac{dA_t}{dt}u_1  
	=\mathcal{L}^{\frac{\alpha}{2}}\frac{W_t-W_{-t}}{2i}u_0+\frac{W_t+W_{-t}}{2}u_1.
	\end{equation*}
By \eqref{LPLP}, we know
	\begin{equation*}
	\|\Delta_jv\|_{L^{r_1}(H^n)}=\left\|\mathcal{L}^{-\frac{\alpha}{2}}\mathcal{L}^{\frac{\alpha}{2}}\Delta_jv\right\|_{L^{r_1}(H^n)}\leq C2^{-j\alpha}\|\mathcal{L}^{\frac{\alpha}{2}}\Delta_jv\|_{L^{r_1}(H^n)}.
	\end{equation*}
Thus we have
	\begin{equation*}
	\|v\|_{L^{p_1}(\mathbb{R},\dot{B}_{r_1,2}^{\rho_1}(H^n))}\leq C\|\mathcal{L}^{\frac{\alpha}{2}}v\|_{L^{p_1}(\mathbb{R},\dot{B}_{r_1,2}^{\rho_1-\alpha}(H^n))}.
	\end{equation*}
Then
	\begin{equation*}
	\begin{aligned}
	&\|v\|_{L^{p_1}(\mathbb{R},\dot{B}_{r_1,2}^{\rho_1}(H^n))}+\|\partial_tv\|_{L^{p_1}(\mathbb{R},\dot{B}_{r_1,2}^{\rho_1-\alpha}(H^n))} \\
	&\leq 
	C\Big(\left\|W_t\mathcal{L}^{\frac{\alpha}{2}}u_0\right\|_{L^{p_1}(\mathbb{R},\dot{B}_{r_1,2}^{\rho_1-\alpha}(H^n))}
	+\left\|W_{-t}\mathcal{L}^{\frac{\alpha}{2}}u_0\right\|_{L^{p_1}(\mathbb{R},\dot{B}_{r_1,2}^{\rho_1-\alpha}(H^n))}  \\
	&+\|W_tu_1\|_{L^{p_1}(\mathbb{R},\dot{B}_{r_1,2}^{\rho_1-\alpha}(H^n))}
	+\|W_{-t}u_1\|_{L^{p_1}(\mathbb{R},\dot{B}_{r_1,2}^{\rho_1-\alpha}(H^n))}\Big).
	\end{aligned}
	\end{equation*}
According to \eqref{WEresult1}, we get
	\begin{equation*}
	\|W_t\mathcal{L}^{\frac{\alpha}{2}}u_0\|_{L^{p_1}(\mathbb{R},\dot{B}_{r_1,2}^{\rho_1-\alpha}(H^n))} \leq \|\mathcal{L}^{\frac{\alpha}{2}}u_0\|_{L^2}\leq\|u_0\|_{\dot{H}^\alpha(H^n)},
	\end{equation*}
	\begin{equation*}
	\|W_tu_1\|_{L^{p_1}(\mathbb{R},\dot{B}_{r_1,2}^{\rho_1-\alpha}(H^n))} \leq \|u_1\|_{L^2(H^n)}.
	\end{equation*}
The same discussion can apply to $W_{-t}$. Now we obtain \eqref{WResult1}.

It remains to prove \eqref{WResult2}. We write
	\begin{equation*}
	w=\int_{0}^{t}A_{t-\tau}f(\tau)d\tau=\int_{0}^{t}\mathcal{L}^{-\frac{\alpha}{2}}\frac{W_{t-\tau}-W_{\tau-t}}{2i}f(\tau)d\tau,
	\end{equation*}
Similarly, we have
	\begin{equation*}
	\partial_tw=\int_{0}^{t}\dfrac{dA_{t-\tau}}{dt}f(\tau)d\tau=\int_{0}^{t}\frac{W_{t-\tau}-W_{\tau-t}}{2}f(\tau)d\tau.	
	\end{equation*}
Therefore, we obtain
	\begin{equation*}
	\begin{aligned}
	&\|w\|_{L^{p_1}([0,T],\dot{B}_{r_1,2}^{\rho_1}(H^n))} \\
	&\leq
	C\left(\left\|\mathcal{L}^{-\frac{\alpha}{2}}\int_{0}^{t}W_{t-\tau}f(\tau)d\tau\right\|_{L^{p_1}([0,T],\dot{B}_{r_1,2}^{\rho_1}(H^n))}\right.\\
	&\qquad\qquad\qquad\qquad\left.+\left\|\mathcal{L}^{-\frac{\alpha}{2}}\int_{0}^{t}W_{\tau-t}f(\tau)d\tau\right\|_{L^{p_1}([0,T],\dot{B}_{r_1,2}^{\rho_1}(H^n))}\right) \\
	& \leq
	C\left(\left\|\int_{0}^{t}W_{t-\tau}f(\tau)d\tau\right\|_{L^{p_1}([0,T],\dot{B}_{r_1,2}^{\rho_1-\alpha}(H^n))}\right.\\
	&\qquad\qquad\qquad\qquad\left.+\left\|\int_{0}^{t}W_{\tau-t}f(\tau)d\tau\right\|_{L^{p_1}([0,T],\dot{B}_{r_1,2}^{\rho_1-\alpha}(H^n))}\right),
	\end{aligned}
	\end{equation*}
and
	\begin{equation*}
	\begin{aligned}
	&\|\partial_tw\|_{L^{p_1}([0,T],\dot{B}_{r_1,2}^{\rho_1-\alpha}(H^n))}  \\
	&\leq C\Big(\left\|\int_{0}^{t}W_{t-\tau}f(\tau)d\tau\right\|_{L^{p_1}([0,T],\dot{B}_{r_1,2}^{\rho_1-\alpha}(H^n))}+\left\|\int_{0}^{t}W_{\tau-t}f(\tau)d\tau\right\|_{L^{p_1}([0,T],\dot{B}_{r_1,2}^{\rho_1-\alpha}(H^n))}\Big).
	\end{aligned}
	\end{equation*}
According to \eqref{WEresult2}, then
	\begin{equation*}
	\left\|\int_0^tW_{t-\tau}f(\tau)\,d\tau\right\|_{L^{p_1}([0,T], \dot{B}^{\rho_1-\alpha}_{r_1,2}(H^n))}\leq C\|f\|_{L^{\bar{p}_2}([0,T],\dot{B}_{\bar{r}_2,2}^{-\rho_2}(H^n))}.
	\end{equation*}
Thus \eqref{WResult2} holds. \\
\qed

As Corollary \ref{SLebesgue}, we obtain the Strichartz inequalities on the Lebesgue space by the embedding relation between the homogeneous Besov space and Lebesgue space on the Heisenberg group.
\begin{corollary}
Let $0<\alpha<2$ and $u$ be the solution of the fractional wave equation \eqref{FWaveEqu}. Suppose $p\in[2N-2\alpha,+\infty)$ and $q$ satisfy
	\begin{equation*}
	\dfrac{\alpha}{p}+\dfrac{N}{q}=\dfrac{N}{2}-\alpha.
	\end{equation*}
Then the following estimate holds true
	\begin{equation*}
	\|u\|_{L^p([0,T],L^q(H^n))} \leq C\left(\|u_0\|_{\dot{H}^\alpha(H^n)}+\|u_1\|_{L^2}+\|f\|_{L^1([0,T],L^2(H^n))}\right),
	\end{equation*}
where the constant $C>0$ does not depend on $u_0$, $u_1$, $f$ or $T$.
\end{corollary}
\begin{remark}
When $\alpha=1$, it corresponds to the results in \cite{BGX2000}. Therefore, our work is more general.
\end{remark}

\noindent\\[2mm]
\subsection{The fourth-order Schr\"{o}dinger Equation}
Finally, we consider the fourth-order Schr\"{o}dinger equation
	\begin{equation}\label{SSchrEqu}
	\begin{cases}
	i\partial_tu+\mathcal{L}^2u+\mathcal{L}u=f,\\
	u|_{t=0}=u_0.
	\end{cases}
	\end{equation}
The solution is given by
	\begin{equation*}\label{solution3}
	u(t)=U_tu_0-i\int_0^tU_{t-\tau}f(\tau)\,d\tau,
	\end{equation*}
where $U_t=e^{it(\mathcal{L}^2+\mathcal{L})}$ and it corresponds to the case when $\phi(r)=r^2+r$. By a simple calculation,
	\begin{equation*}
	\phi'(\lambda)=2r+1, \quad
	\phi''(\lambda)=2.
	\end{equation*}
We know $\phi$ satisfies (H1)?(H4) with $m_1=\alpha_1=\alpha_2=2$, $m_2=1$. According to Theorem\ref{ResultTime}, we get
	\begin{equation*}
	\|U_t \varphi_j\|_{L^\infty (H^n)} \leq C_n|t|^{-\frac{1}{2}}2^{j(N-2)}.
	\end{equation*}

\begin{theorem}
If $u$ is the solution of the free fourth-order Schr\"{o}dinger equation \eqref{FSchrEqu} (i.e. with $f=0$), then there exists a constant $C>0$ which does not depend on $u_0$ or $t$ such that
	\begin{equation*}
	\|u(t)\|_{L^\infty (H^n)}\leq C|t|^{-1/2}\|u_0\|_{\dot{B}^{N-2}_{1,1}},
	\end{equation*}
and the decay is sharp in time.
\end{theorem}

Since the proof of the inequality is similar to the case of the fractional Schr\"{o}dinger equation, we only show the sharpness of the decay. For the same $Q$, $D_0$ and $u_0$, and if $u$ is the solution of \eqref{SSchrEqu} when $f=0$, we have
	\begin{equation*}	
	u(z,s,t)=e^{it(\mathcal{L}^2+\mathcal{L})}u_0=C\int_{D_0}e^{-i\lambda s}e^{it((4n\lambda)^2+4n\lambda)}Q(\lambda)e^{-\lambda|z|^2}\lambda^n \,d\lambda.
	\end{equation*}
In particular,
	\begin{equation*}
	u(0,ts,t)=C\int_{D_0}e^{it(-\lambda s+(4n\lambda)^2+4n\lambda}Q(\lambda)\lambda^n \,d\lambda.
	\end{equation*}
This oscillating integral has a phase $\Psi(\lambda)=-\lambda s+(4n\lambda)^2+4n\lambda$ with a unique critical point $\lambda_0=\dfrac{s-4n}{32n^2}$, which is not degenerate. Applying Lemma \ref{asymptotic}, for $s_0=32n^2+4n$, we have
	\begin{equation*}
	u(0,ts_0,t) \sim C|t|^{-1/2}.
	\end{equation*}

\begin{proposition}\label{SResultSB}
 For $i=1,2$, let $p_i, r_i\in[2,\infty]$ and $\rho_i\in\mathbb{R}$ such that
	\begin{equation*}
	\begin{aligned}
	&(i)\,\dfrac{2}{p_i}+\dfrac{1}{r_i} \leq \dfrac{1}{2},   \\
	&(ii)\,\frac{N}{r_i}+\frac{4}{p_i}-\rho_i=\frac{N}{2},
	\end{aligned}
	\end{equation*}
Then the following estimates are satisfied:
	\begin{equation*}
	\|U_tu_0\|_{L^{p_1}(\mathbb{R}, \dot{B}_{r_1,2}^{\rho_1}(H^n))}\leq C\|u_0\|_{L^2(H^n)},
	\end{equation*}
	\begin{equation*}
	\left\|\int_0^tU_{t-\tau}f(\tau)\,d\tau\right\|_{L^{p_1}([0,T], \dot{B}^{\rho_1}_{r_1,2}(H^n))}\leq C\|f\|_{L^{\bar{p}_2}([0,T],\dot{B}_{\bar{r}_2,2}^{-\rho_2}(H^n))},
	\end{equation*}
	where the constant $C>0$ does not depend on $u_0$, $f$ or $T$.
\end{proposition}

\begin{theorem}\label{Fsolution}
Under the same hypotheses as in Theorem \ref{SResultSB}, the solution of
the fourth-order Schr\"{o}dinger equation \eqref{SSchrEqu} $u$ satisfies the following estimate
	\begin{equation*}
	\|u\|_{L^{p_1}([0,T], \dot{B}^{\rho_1}_{r_1,2}(H^n))} \leq C\left(\|u_0\|_{L^2(H^n)}+\|f\|_{L^{\bar{p}_2}([0,T],\dot{B}_{\bar{r}_2,2}^{-\rho_2}(H^n))}\right),
	\end{equation*}
where the constant $C>0$ does not depend on $u_0$, $f$ or $T$.
\end{theorem}

\begin{corollary}\label{FLebesgue}Suppose $p\in[2N-4,+\infty]$ and $q$ satisfy
	\begin{equation*}
	\dfrac{4}{p}+\dfrac{N}{q}=\dfrac{N}{2}-1?
	\end{equation*}
If $u$ is a solution of the fourth-order Schr\"{o}dinger equation \eqref{SSchrEqu}, then the following inequality holds true
	\begin{equation*}
	\|u\|_{L^p([0,T],L^q(H^n))} \leq C\left(\|u_0\|_{\dot{H}^1(H^n)}+\|f\|_{L^1([0,T],\dot{H}_1(H^n))}\right)?
	\end{equation*}
where the constant $C>0$ does not depend on $u_0$, $f$ or $T$.
\end{corollary}

The proofs of Proposition \ref{SResultSB}, Theorem \ref{Fsolution}, Corollary \ref{FLebesgue} are similar to those in the case of the fractional Schr\"{o}dinger case. We omit it here.

\noindent\\[4mm]
\noindent\bf{\footnotesize Competing interests}\quad\rm
{\footnotesize The author declares to have no competing interests.}\\[4mm]
\noindent\bf{\footnotesize Acknowledgements}\quad\rm
{\footnotesize The work is supported by the National Natural
Science Foundation of China (Grant No. 11701452),  the China Postdoctoral Science Foundation (Grant No. 2017M623230) and the Natural Science Foundation of Shaanxi Province (Grant No. 2020JQ-112).}\\[4mm]

\noindent{\bbb{References}}

\end{document}